\documentclass[12pt]{article}

\usepackage{amsthm,amsmath,amssymb,xy}
\xyoption{all}

    \newtheorem{theorem}{Theorem}[section]
    \newtheorem{proposition}[theorem]{Proposition}
    \newtheorem{lemma}[theorem]{Lemma}
    \newtheorem{corollary}[theorem]{Corollary}

    \newtheorem{nul}[theorem]{}
\DeclareMathOperator{\Cliff}{Cliff}
\DeclareMathOperator{\Hom}{Hom}
\DeclareMathOperator{\Ext}{Ext}

\DeclareMathOperator{\Imag}{Im}
\DeclareMathOperator{\Ker}{Ker}

\def\O{{\mathcal{O}}}
\def\C{{\mathbb {C}}}
\def\P #1{{\mathbb P}^{#1}}
\def\G{{\mathbb {G}}}
\def\F{{\mathbb {F}}}
\def\S{{\mathcal S}_+}
\def\qed{\hfill $\blacksquare$}

\numberwithin{equation}{theorem}

\title{\bf Extensions of Instantons}
\author{
\begin{tabular}{l} Cristian Anghel
\end{tabular}
and
\begin{tabular}{l} Nicolae Manolache
\end{tabular}}
\date{}

\begin{document}
\maketitle
\begin{abstract}
In this paper one studies a correspondence between extensions $E$ of instanton
bundles $F$ on $\P {3}$ (by line bundles or by twisted bundles of forms) and
potentials (connections) on certain bundles on the grassmanian ${\mathbf Gr}(2,4)$,
or the sphere $S^4$. Of special interest is the case of extensions by $\O (-1)$,
when the $E$'s are stable. One shows also how to obtain the monads of the $E$'s,
knowing those of the $F$'s. The main ingredients are the Atiyah-Ward correspondence
(cf. \cite{AW}) and a theorem of Hitchin (cf. \cite{H}).
\end{abstract}

\section{Introduction.}\label{intr}
\noindent

   In the fundamental papers \cite{W}, \cite{AW} of Ward and Atiyah and Ward
   a 1-1 correspondence between
   {\em \{vector bundles on the grassmanian ${\mathbf Gr}(2,4)$ (respectively $S^4$)
    with an antiself-dual connection\}}
   and {\em \{holomorphic vector bundles on $\P{3}$
   (respectively vector bundles on $\P{3}$ with a positive real form)\}}
   is obtained.

   In particular, the bundles on $\P{3}$ in this correspondence have  necessarily
   $c_{1}=0$.
   A natural question, which in fact the second named author learned from
   a discussion with  M. F. Atiyah in 1981 in Bucharest, is to find a similar correspondence
   between vector bundles on $\P{3}$ with arbitrary $c_1$ and some objects of
   differential geometry on ${\mathbf {Gr}}(2,4)$ or $S^4$.
   In this direction, an answer was given by Leiterer in \cite{L},
   who established a correspondence between arbitrary holomorphic
   vector bundles on $\P{3}$ and holomorphic bundles over
   ${\mathbf {Gr}}(2,4)$ with a certain differential operator which he named
   {\em partial connection}.

   The aim of this paper is to describe a construction by which to
   certain holomorphic vector bundles $E$ on $\P{3}$ with
   any $c_{1}$, namely extensions of an instanton bundle, one makes to
   correspond  bundles on ${\mathbf {Gr}}(2,4)$ equiped with certain {\em connections}
   constructed from data associated with the instanton. The ingredients of the
   construction are the Atiyah-Ward correspondence  and the choice of a spinor on
   ${\mathbf {Gr}}(2,4)$ coupled with the Yang-Mills potential obtined by the twistor
   transform. Besides these, we use an elementary construction
   of a connection on a direct sum of vector bundles when there are known
   two connections on the two bundles  and a section on their tensor
   product.

   When $c_{1}(E)=-1$ one uses a theorem of Hitchin (cf. \cite{H}) to obtain
   a bijection between certain vector bundles on $\P{3}$ and a special class of connections
   in an associated vector bundle over $S^4$.

   At the end of {\bf Section \ref{ConnConstr}}  we extend our method for obtaining
   families of connections on bundles over ${\mathbf {Gr}}(2,4)$ associated to families
   of pairs of a cohomology class of an instanton bundle on
   $\P{3}$ and a spinor field on ${\mathbf {Gr}}(2,4)$ coupled with the
   Yang-Mills potential which corresponds to the instanton.

   The Atiyah-Ward correspondence makes the construction of instantons on $\P{3}$
   of a special relevance. A basic tool is the {\em monads method}
   of G. Horrocks (cf. \cite{OSS}, \cite{A} ).
   In {\bf Section \ref{Mon}} it is shown how to obtain monads of our vector bundles $E$
   from the monads of the instantons $F$.

   It should be pointed that, in "physicists" terms, our construction
   can be summarised  by saying that, from a Yang-Mills potential and
   a spinor field with arbitrary helicity, one constructs a potential
   (i.e. connection) on a vector bundle of higher rank.

{\bf Notation}
We explain here some of the notation we use:

\begin{description}
  \item{$\bullet$} $\P{} :=\P{} (U^\ast)\cong \P{3}=$ projective space of lines in $U$,
   where  $U$ is a $\C$-vector space of dimension $4$
  \item{$\bullet$} $\G:={\mathbf {Gr}}(2,4)=$ the Grasmannian of  $2$-dimensional subspaces
       in a vector space of dimension $4$
  \item{$\bullet$} $\S :=$ positive spinor bundle on $\G$, i.e. universal subbundle on
  $\G$
  \item{$\bullet$} $\F :=$ incidence variety $\{(x,l) : x\in \P{3}, l\in \G
        \hbox{ with the property }
        x\in L_l\}$, $L_l$ being the line in $\P{}$ represented by $l$
  \item{$\bullet$} $S^4 :=$ the $4$-sphere
  \item{$\bullet$}  $F$ a mathematical instanton bundle on $\P {3}$, i.e a
        self dual vector bundle on $\P{}$  of arbitrary rank, such that $H^1(F(-2))=0$
  \item{$\bullet$}  $VE:= V\otimes E$, where $V$, $E$ are vector bundles, even
       when $V$ is a vector space interpretted as a constant vector bundle
  \item{$\bullet$} $V\O$  for $V$ a vector space considered as constant vector bundle
\end{description}

\section{Construction of the connection}\label{ConnConstr}

\subsection{General construction}\label{GenConstr}

  Let $G$ be a differentiable manifold and $F_{1}$, $F_{2}$ vector
  bundles over $G$ endowed with connections ${\nabla}_{1}$,
  ${\nabla}_{2}$. Let ${\varphi}_{1} \in \Gamma (F_{1}\otimes
  {F_{2}}^{\ast})$ and ${\varphi}_{2} \in \Gamma (F_{2}\otimes
  {F_{1}}^{\ast})$ be two differentiables sections, where by $\ast$
  we denote the dual of the corresponding vector bundle. Let $F =
  F_{1}\oplus F_{2}$. To the quadruple
  $({\varphi}_{1},{\varphi}_{2},{\nabla}_{1},{\nabla}_{2})$ we associate
  the following connection $\nabla$ on $F$: relative to the
  splitting $F = F_{1}\oplus F_{2}$, $\nabla$ has the block
  decomposition
\[
\begin{pmatrix} \nabla _1 & \nabla _{12}\varphi _1   \\
                           &                          \\
                \nabla_{21}\varphi _2  & \nabla _2
\end{pmatrix}
\]

  where ${\nabla}_{12}$ and ${\nabla}_{21}$ are the canonical
  induced connections on $F_{1}\otimes {F_{2}}^{\ast}$ and
  $F_{2}\otimes {F_{1}}^{\ast}$. If we change the splitting it is
  obvious that the connection is changed by a gauge transformation
  and so we obtain the following simple lemma:
\begin{lemma}\label{genConstr}
  In the situation from above to any quadruple
  $({\varphi}_{1},{\varphi}_{2},{\nabla}_{1},{\nabla}_{2})$ it is
  associated a unique gauge equivalence class of connections on
  $F$.

  \qed
\end{lemma}

\subsection{Modification of the Atiyah-Ward transform}

  Recall the notation sumarized in {\bf Introduction}:
  $\G$ will denote the grassmanian of $2$-planes in
  $\C^4$, $\P{}$ the projective $3$-space and
  $\F$ the (flag) incidence variety. We shall use the
  following twistor diagram (cf., for instance, \cite{A}, \cite{WW}):

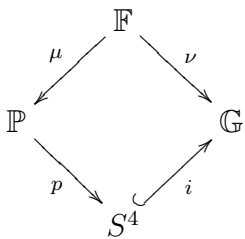
\[
\xymatrix{
           & \F\ar[dr]^\nu \ar[dl]_\mu & \\
           \P{}\ar[dr]_p & & \G   \\
           & S^4\ar@{^{(}->}[ur]_i &
}
\]

  Recall that $\G$ is realized as a quadric $Q_4$ in $\P{5}$
  via the Pl\"{u}cker embedding (famous Klein representation of all  lines in $\P{3}$).
  A real structure on $\P{3}$ induces a real structure on $Q_4$.
  Taking the real structure defined by
\[
\sigma : (z_1,z_2,z_3,z_4)\mapsto (-\bar{z}_2,\bar{z}_1,-\bar{z}_4,\bar{z}_3)
\]
  one shows that although the above map has no fixed points in $\P{}$,
  it has fixed lines, which are called {\em real lines}. The real lines are
  parametrized by  $S^4\hookrightarrow \G$, and $\G$ is the complexification
  of $S^4$. The map $p$ associates to any point $x$ in $\P{}$ the
  point in $S^4 \subset \G$ representing the real line determined by $x$ and
  $\sigma (x)$.

  Let $F$ be a rank-$2$ instanton bundle on $\P{}$ i.e. a
  mathematical instanton trivial on all the
  real lines. Let $ \mathcal{U} \subset G$ be the open subset which
  corresponds to all lines in $\P{}$ on which $F$ is trivial. $
  \mathcal{U}$ will therefore necessarily contain $S^{4}$. Let
  ${\mathcal{U}}'= \mu \circ {\nu}^{-1}(\mathcal{U})$.

  It is well known that $H^1(F(k))=0$ for $k\le -2$ and thus we
  start the construction with an element $e \in
  H^1(F(k))\cong \Ext (F,\O (k))$ where $k \ge -1$. To $e$  it corresponds an
  extension of the form:
\[
  0 \to \O(k)  \to E \to F \to 0
\]

  and we want to apply the twistor transform to the preceding
  sequence. First of all, the condition $k \ge -1$ implies that the
  restriction of $E$ to all real lines is of the form $\O(k)
  \oplus \O \oplus \O$, because on $\P{1}$ we
  have $\Ext^{1}( \O \oplus \O, \O(k)) =0$ if $k
  \ge -1$. So, applying $ {\mu} ^{\ast}$ and $ {\nu} _{\ast}$ one
  obtains the following sequence on an open set $ \mathcal{U}$
  containing $S^{4}$ which is exact also on the right because on
  $\P{1}$, $ H^1(\O(k))=0$ if  $k \ge -1$:

\[
 0 \to {\mathcal S}_{k}  \to \tilde{E} \to \tilde{F} \to 0.
\]

  Now, ${\mathcal S}_{k}$ is the $k$-symmetric power of the positive
  spinor bundle over $\G$ and it has therefore the canonical
  Levi-Civitta connection denoted by ${\nabla}_{k}$. $\tilde{F}$ is
  the Atiyah-Ward transform of $F$ and it has the Yang-Mills
  connection $\tilde {\nabla}$. According to {\bf \ref{GenConstr}}, for
  constructing a connection on $\tilde{E}$, we need also two sections

\[
  {\varphi}_{1} \in \Gamma ({\tilde{F}} ^{\ast} \otimes {\mathcal S}_{k})
  \hbox{ and } {\varphi}_{2} \in \Gamma ({\tilde{F}}\otimes {{\cal
  S}_{k}}^{\ast} ).
\]
  But, using the canonical metric on ${\mathcal S}_{k}$ and the fact
  that ${\tilde{F}} \cong {\tilde{F}} ^{\ast}$, any section
  ${\varphi}_{1}$ will automatically determine a unique section
  ${\varphi}_{2}= {{\varphi}_{1}} ^{\ast}$. The preceddings remarks and
  the {\bf Lemma} \ref{genConstr} imply therefore the following result:

\begin{theorem}\label{modifAW}
  Consider the following set of data: $F$ a rank 2 instanton bundle
  on $\P{3}$, $e \in H^1(F(k))$ where $k \ge -1$,
  ${\varphi} \in \Gamma ({\tilde{F}} ^{\ast} \otimes {\mathcal S}_{k} )$.
  Then on $\tilde{E}$, the twistor transform of the extension of $F$
  corresponding to $e$, there exist a canonical connection
  associated to the triple $(F, e,{\varphi})$.

  \qed
\end{theorem}

\subsection{The case $k = -1$}\label{-1}

  The value $k = -1$ is special for at least two reasons:

  First of all, in this case the $E$'s are stable, because $H^0(E)=H^0(E^\ast (-1))=0$
  (cf. \cite{OSS} Remark 1.2.6 Chap. II). This observation suggests that an
  interesting problem would be the study of the geometry of these families of
  extensions of instantons in the ambient moduly spces.

  Secondly, as it was proved by Hitchin
  in \cite{H}, $H^1(F(-1))$ can be identified with the subspace of sections in
  $\Gamma (\S\otimes \tilde{F})$ which satisfy the Dirac equation
  coupled with the Yang-Mills potential on $\tilde{F}$, where now,
  $\tilde{F}$ denotes the restriction of the twistor transform of
  $F$ to $S^{4}$, and $\S$ the positive spinor bundle over
  $S^{4}$.

  But now one can apply directly to $\tilde{F} \oplus \S$
  the construction from {\bf Section \ref{GenConstr}}, whitout using  $E$ directly,
  but using only its extension class $e \in H^1(F(-1))$ for
  obtaining the section $\psi _2  \in \Gamma (\S \otimes
  \tilde{F})$ by \cite{H}.

  Conversely, let $\mathcal{A}$ be the space of connections on
  ${\tilde{F}} \oplus \S$ such that in the block
  decomposition
\[
  \begin{pmatrix} \widetilde \nabla  & \Theta _1   \\
                           &                          \\
                \Theta_2  & \nabla _1
\end{pmatrix}
\]
$\Theta _ 1$ and $\Theta _2$ have the followings properties:
\[
\begin{array}{rlcl}
    (i)  & \Theta _1             &=&  \nabla (\psi _1)  \\
   (ii)  & \Theta _ 2            &=&  \nabla (\psi _ 2)  \\
   (iii) & \psi _ 2              &=&  \psi _1^\ast   \\
   (iv)  & \Cliff(\Theta _ 1)    &=&  0
\end{array}
\]

  where $\Cliff$ denotes the Clifford multiplication.

  For any connection in $\mathcal{A}$ one obtains a stable bundle on
  $\P{3}$ and so we have the following:

\begin{corollary}
  There exists a canonical bijection between the class of stable vector
  bundles on $\P{3}$  which are extensions of the
  instanton bundle $F$ by $\O(-1)$ and connections in the
  class $\mathcal{A}$.

  \qed
\end{corollary}

\subsection{Extensions by twisted holomorphic forms}\label{ExtForms}

  The aim of this section is to find an analogue of the above
  {\bf Corollary} when one starts with an extension of the following form:
\[
  0 \to \Omega^1  \to E \to F \to 0.
\]
  Let $e \in H^1( \Omega^1 \otimes F)$ be the class of this
  extension. Using the standard sequence on $\P{}$,
\bigskip

\makebox[\textwidth][s]{
\hfill
$ 0 \rightarrow  \Omega^1  \rightarrow U \O(-1) \rightarrow \O \rightarrow 0 $
\hfill $(\ast)$
}
\bigskip

 (we recall that $U = H^0(\O(1))$), and tensoring it by $F$,
  one obtains the following part of the cohomology sequence:
\[
  H^0(F) \to H^1( \Omega^1 F) \to U H^1(F(-1)) \to H^1(F).
\]
  The first term on the left is $0$ by the stability of the
  instanton bundle $F$. Denote by $m$ the last arrow on the
  right, which in fact is a component of the structural multiplication
  of the $\oplus H^0(\O (n))$-module $\oplus H^1(F(n))$.
  The cohomology class $e$ will determine a quadruple of
  elements $(e_i)$ so that $m ((e_i))=0$ and conversely
  any such quadruple will determine a unique extension of $F$ by
  $\Omega ^1$. Moreover, using \cite{H} as in the preceeding
  section, any such quadruple  corresponds to four Dirac spinor
  on $S^4$ coupled with the Yang-Mills potential $\tilde{F}$.
  Therefore we obtain a canonical connection on ${\tilde{F}} \oplus
  U \S$ in which the  off diagonal terms
  ${\Theta}_{1}$ and ${\Theta}_{2}$ have the followings properties:

\[
\begin{array}{rlcl}
    (i)   & \Theta _1            &=& \nabla (\psi _1) \\
   (ii)   & \Theta _2            &=& \nabla (\psi _ 2)  \\
   (iii)  & \psi _ 2             &=& \psi _ 1^\ast     \\
   (iv)   & \Cliff( \Theta _ 1 ) &=&  0 \\
   (v)    & m(\psi _1)      &=&  0
\end{array}
\]

  where $\Cliff$ denotes the Clifford multiplication, $\psi _1$ is
  a quadruple of Dirac spinors and we identify according to  \cite{H} the
  space of Dirac spinors with the corresponding cohomology group on
  $\P{}$. Conversely, for any connection with these
  properties we obtain a unique vector bundle on $\P{}$
  which is an extension of the instanton $F$ by  $\Omega ^1$
  and
  so we have the following:

\begin{corollary}
  There exists a canonical bijection between the classes of extensions
  on $\P{3}$ of instantons by ${\Omega}^{1}$ and
  connections on $S^{4}$ which satisfy the five properties above.

\qed
\end{corollary}

  We consider now extensions of the form:

\[
0 \to {\Omega}^{2}(1)  \to E \to F \to 0     \ .
\]
  We shall use the following standard exact sequence:

\[
  0 \to \Omega ^2  \to \wedge ^2 U \O(-2)
  \to \Omega ^1 \to 0.
\]

  By tensoring the last
  sequence with $F(1)$ one obtains the following exact cohomology sequence:

\[
  H^0(F \Omega ^1(1)) \to H^1( \Omega ^2
  F(1)) \to \wedge ^2 U  H^1(F(-1)) \to H^1(
  \Omega ^1  F(1)).
\]

  The first term is $0$, as one sees multiplying the exact sequence $(\ast)$
  with $F(1)$, taking the cohomology and using the stability of $F$.
  The last arrow on the right is the multiplication
  $m: H^0(\Omega^1 (2))\otimes H^1(F(-1))\rightarrow H^1(\Omega ^1F(1))$.

  Therefore one obtains the following analogue of the preceeding {\bf Corollary}:

\begin{corollary}
  There exist a canonical bijection between the class of extensions
  on $\P{3}$ of instantons by $\Omega^2(1)$ and
  connections of  $\tilde{F} \oplus {\wedge}^{2}U \otimes \S$ over
  $S^4$ which satisfy the analogues of the five properties above.

\qed
\end{corollary}

\begin{nul}{\bf Remark. } Let now $F$ be a rank two holomorphic
  bundle over $\P{3}$ which {\em is not} an instanton, namely, it satisfies
  the following:
\[
 \hspace{1cm} F \hbox{ is trivial on a line and } H^1(F(k))
 \neq 0 \hbox{ for some } k \le
  -2.
\]
  For such a $F$ one can  still consider extensions of the form:
\[
  0 \to \O(k)  \to E \to F \to 0,
\]
  One restricts  $E$ to the set ${\mathcal U}'=
  \mu(\nu^{-1}(\mathcal U ))$, where $\mathcal U$ is the open set
  in $\G$ which parametrise the lines on which $F$ is trivial, and
  one takes the  cohomology class $e \in H^1(\mathcal U, F(k))$.
  This cohomology group can be
  identified by \cite{WW} with the space of massless fields on
  $\mathcal U $ of helicity $=\frac{-k-2}{2}$, coupled with the
  Yang-Mills potential obtined on $\tilde{F}$ by the Penrose
  transform. As in the precedings sections it will be a canonical
  connection on ${\mathcal S}_{-k-2}\oplus \tilde{F}$.
\end{nul}

\subsection{Towers of connections}

  The aim of this section is to generalize the construction in $2.2$
  in order to obtain towers of connections on the grassmanian.

  We start with the following set of data, for $k\ge -1$: $F$ a rank-2 instanton
  bundle on $\P{3}$, $e_{k-i}\in H^1(F(i))$ for all $i$, $-1\le i \le k$ and
  a set of sections
  ${\varphi}_{k-i} \in \Gamma ({\tilde{F}} ^{\ast} \otimes {\mathcal
  S}_i )$. We describe only the first steps of the construction
  and skip the inductive argument which is obvious.

  For $i=k$ the construction is the same as in $2.2$ with $e$
  raplaced by $e_0$, ${\varphi}$ by $\varphi _0$ and $E$ by
  $E_0$. We denote by $\nabla _0$ the corresponding connection
  on $\tilde{E}_0$.

  Let be $i=k-1$. Dualizing and  twisting the exact sequence
\[
  0 \to \O(k)  \to E_0 \to F \to 0,
\]

  one gets the exact sequence:
\[
  0 \to F^\ast(k-1)  \to E_0^{\ast}(k-1) \to \O(-1)
\to 0.
\]
  By taking the cohomology one obtains:

\[
  H^1(F^{\ast}(k-1)) \cong H^1(E_0^\ast(k-1))
\]
  and, as $F\cong F^\ast$, we find that $e_{k-1}$
  determines an extension of the following type:
\[
  0 \to \O(k-1)  \to E_1 \to E_0 \to 0.
\]
  On the grassmanian we will have therefore the following sequence:
\[
0 \to {\mathcal S}_{k-1}  \to \tilde{E}_1 \to \tilde{E}_0 \to 0.
\]
  To construct a connection on $\tilde{E}_1$ one needs a section
  in $\Gamma (\tilde{E}_0^\ast \otimes {\mathcal S}_{k-1} )$.
  To obtain it, we use the element ${\varphi}_{k-1} \in \Gamma
  (\tilde{F}^\ast \otimes {\mathcal S}_{k-1} )$ and the canonical
  homomorphism
\[
  \Gamma (\tilde{F} ^\ast \otimes \mathcal S _{k-1} ) \to \Gamma
  (\tilde{E}_0 ^\ast \otimes \mathcal S_{k-1} ).
\]
  We denote by $\tilde{\varphi}_{k-1}$ the image of $\varphi _{k-1}$ by the preceeding
  homomorphism. Following the construction in
  {\bf Section \ref{ConnConstr}} we obtain a canonical connection on $\tilde{E}_1$.
  This construction can be performed inductively and gives:
\begin{corollary}
  For $k\ge -1$, to the data: $F$ a rank-2 instanton
  bundle on $\P{3}$,  $e_{k-i}\in
  H^1(F(i))$ for $-1\le i \le k$ and sections
  $\varphi_{k-i} \in \Gamma (\tilde{F} ^\ast \otimes {\mathcal
  S}_i )$ there corresponds a canonical tower of connections on
  the Penrose transformed of the successive extensions of $F$.

\qed
\end{corollary}
\begin{nul}
{\bf Remark.} Two interesting special cases arise:

1) the preceeding construction can be modified to extend $F$ successively by $\O$

2) extend $F$ succesively by $\O(k_i)$, with incresing $k_i$'s.

\end{nul}
\section{Monads for Extensions of Instanton Bundles on $\P {3}$}\label{Mon}
\begin{nul}
     Vector bundles can be constructed using the Horrock's monads, i.e.
     complexes of vector bundles of the type
\[
  \xymatrix{
  0\ar[r] & \mathcal{W}_1\ar[r]^A & \mathcal{V} \ar[r]^B & \mathcal{W}_2\ar[r] & 0  \  ,
  }
\]
    where $A$ and $B$ are injective, respectively surjective, homomorphisms of
    vector bundles. The object of homology $\Ker B /\Imag A $ is a vector bundle.

    Acording to a theorem of Barth (cf. \cite{BH}, or \cite{OSS}), a symplectic or
    orthogonal vector bundle $F$ on $\P {3}$, of rank $r$ which  satisfies the conditions:

\begin{description}
  \item{$\bullet$} are trivial on some line $\ell$ in  $\P {3}$,
  \item{$\bullet$} $H^1(F(-2))=0$ \ ,
\end{description}

    conditions fulfilled by the the instantons we considered so far,
    is given by a monad of the special type:
\[
  \xymatrix{
   0\ar[r] & W\O(-1)\ar[r]^A & V\O \ar[r]^{A^\ast } & W^\ast \O (1) \ar[r] & 0
   }
\]
    where $W$, $V$ are vector spaces of dimensionss $n$, respectively $2n+r$ ,
    $A$ is a $(2n+r)\times n$-matrix of linear forms of constant rank $n$ on $\P{3}$, $W^{\ast }$
    is the dual of $W$, $A^\ast $ is the dual of $A$.

    To the monad one associates the so-called {\em display}, i.e. the commutative diagram
    with exact rows and columns:
\[
\xymatrix{
    & &0\ar[d] & 0 \ar[d] &  \\
    0\ar[r] & W\O(-1)\ar[r]\ar@{=}[d] &Q^\ast\ar[r]\ar[d] & F \ar[r]\ar[d] &0 \\
    0\ar[r] & W\O(-1)\ar[r]^A &V\O\ar[r]\ar[d]^{A^\ast } & Q \ar[r]\ar[d] &0  \\
    & & W^\ast \O (1) \ar@{=}[r]\ar[d]& W^\ast\O (1)\ar[d] & \\
    & & 0 & 0
}
\]
    It is interesting to interpret the data in the above monad in terms of the cohomology of
    $F$, (cf. \cite{DM} or \cite{A}):
\begin{eqnarray}
    W^{\ast }\cong H^1(F(-1))   \nonumber \\
    V\cong H^1(F\Omega ^1) \nonumber
\end{eqnarray}
    and the linear map $A^{\ast }$ can be identified with the canonical map
\[
    H^1(F\Omega ^1) \rightarrow H^1(F(-1))\otimes H^0(\O (1))  \ ,
\]
    which comes from the cohomology applied to the Euler exact sequence
\[
    0\rightarrow \Omega ^1 \rightarrow H^0(\O (1))\otimes \O(-1) \rightarrow \O \rightarrow 0
\]
    multiplied (tensorially) by $F$.
\end{nul}
\begin{lemma} \label{GenExt}
    Let $F$ be an instanton bundle on $\P{3}$ and
\[
   \xymatrix{
   0\ar[r] & W\O(-1)\ar[r]^A & V\O \ar[r]^{A^\ast } & W^\ast \O (1) \ar[r] & 0
   }
\]
   its monad.  If $K$ is a vector bundle on $\P{3}$ and $f \in \Hom (W\O (-1),K)$
   is mapped to $e \in H^1(F\otimes K)\cong \Ext ^1(F,K)$ by the canonical homomorphism
\[
  \Hom (W\O (-1),K) \cong H^1(F(-1))\otimes H^0(K) \rightarrow H^1(F\otimes K) \  ,
\]
  then the monad
\medskip

\makebox[\textwidth][s]{\hfill
    $
    \xymatrix{
    0\ar[r] & W\O(-1)\ar[r]^{\begin{pmatrix}A \\ f\end{pmatrix}} & V\O \oplus K
    \ar[rr]^{\begin{pmatrix}A^\ast & 0\end{pmatrix}} && W^\ast \O (1) \ar[r] & 0
    } $
    \hfill $(\ast \ast)$
    }
\medskip

    defines a vector bundle $E$ which is the extension of $F$ by $K$ :
\[
   0 \rightarrow K \rightarrow E \rightarrow F \rightarrow 0  \ ,
\]
    corresponding to the element $e \in \Ext ^1(F,K)$.
\end{lemma}
\begin{proof}
    Exercise in homological algebra.
\end{proof}
\begin{proposition} \label{Ext}
   (i) All vector bundles which are extensions of an instanton by an $\O (k)$
          are given by a monad  of the form:
\[
  \xymatrix{
  0\ar[r] & W\O(-1)\ar[r]^{\begin{pmatrix}A \\ f \end{pmatrix}} & V\O \oplus \O(k)
   \ar[rr]^{\begin{pmatrix} A^\ast  &  0 \end{pmatrix}} && W^\ast \O (1) \ar[r] & 0
  }
\]
where $f$ is a row vector of $\dim W$ forms of degree $k+1$. \medskip

   (ii)  Any extension of an instanton $F$
       by $\Omega ^1(k)$ with $k\ge 1$ or by $\Omega ^2(k)$ with $k\ge 2$ is given by a monad
       of the shape $(\ast \ast )$ . \medskip

\end{proposition}
\begin{proof}
   (i) The map
\[
\Hom (W\O(-1),\O(k))\cong \Hom (\O(k+1))\otimes H^1(F(-1))\rightarrow H^1(F(k))
\]
     is surjective, the $k[X_0,\ldots ,X_3]$-module $\bigoplus _{i\ge -1} H^1(F(i))$
     being generated by $H^1(F(-1))$ (cf. \cite{A} or \cite{OSS}) \bigskip

  (ii)  Multiplying the monad of $F$ by $\Omega = \Omega ^1(k)$, respectively by
        $\Omega = \Omega ^2(k)$
        and taking the cohomology, one obtains the following commutative and exact diagram : \bigskip
\[
\xymatrix{
    H^0(W^\ast \Omega (k+1))\ar@{=}[r]\ar[d] & H^0(W^\ast \Omega (k+1))\ar[d] &  \\
    H^1(Q^\ast \Omega (k))\ar[r]\ar[d] &  H^1(F \Omega (k))\ar[r] &  0   \\
    H^1(V\Omega (k))=0 & &
}
\]
which shows  that the map
\[
 \Hom (W\O (-1),\Omega(k))=H^0(\Omega (k+1))\otimes H^1(F(-1))\rightarrow H^1(F\Omega (k))
\]
is surjective.
\end{proof}
\section*{}  {\bf Acnowledgements.} During the preparation of this paper the first named author
was supported by a DFG grant. He expresses his thanks to DFG. Both authors thank to
Oldenburg University, especially to Professor Udo Vetter, for hospitality.

\bigskip
\begin{tabular}{ll} Cristian Anghel & Nicolae Manolache  \\
Institute of Mathematics   & Institute of Mathematics \\
of the Romanian Academy,   & of the Romanian Academy   \\
P.O. Box 1-764, RO-70700,  & P.O. Box 1-764, RO-70700,  \\
Bucharest, Romania,         & Bucharest, Romania,\\
e-mail: cristian.anghel@imar.ro &  e-mail: nicolae.manolache@imar.ro \\
  &  current address: Fachbereich Mathematik, \\
  & Universit\"{a}t Oldenburg, \\
  & Pf. 2503, D-26111 Oldenburg, Germany, \\
  & e-mail: nicolae.manolache@uni-oldenburg.de
\end{tabular}


\begin{thebibliography}{WWW}

\bibitem[A]{A} Atiyah, M. F.:{\em Geometry of Yang-Mills Fields}, Academia Nazionale dei Lincei,
    Pisa 1979

\bibitem[AW]{AW} Atiyah, M. F. and  Ward, R. S. : {\em
    Instantons and Algebraic Geometry}, Comm. Math. Phys. {\bf 55}, 111-124
    (1977)
\bibitem[BH]{BH} Barth, W. and Hulek, K.: {\em Monads and Moduli of Vector Bundles},
    Manuscripta Math. {\bf 25},323-347 (1978)

\bibitem[DM]{DM} Drinfeld, V. G. and Manin, Yu. I.: {\em Instantons and Sheaves on
     $\C \P{3}$ }, Funk. Analiz., {\bf 13}, 59-74 (1979)

\bibitem[H]{H} Hitchin, N. J. :  {\em Linear Field Equations on Self-dual Spaces},
     Proc. London Math. Soc. {\bf 43}, 133-150 (1980)

\bibitem[Ho1]{Ho1} Horrocks, G.: {\em Vector Bundles on the Punctured Spectrum of a Local Ring},
    Proc. london Math. Soc. 83), {\bf 14}, 689-713 (1964)

\bibitem[Ho2]{Ho2} Horrocks, G.: {\em A Construction of Locally Free Sheaves}, Topology
    {\bf 7} 117-120 (1968)

\bibitem[L]{L} Leiterer, J : {\em The Penrose Transform for Bundles
       Non-trivial on the General Line},
       Math. Nachr. {\bf 112}, 35-67 (1983)

\bibitem[OSS]{OSS} Okonek, C. , Schneider, M. and Spindler, H.: {\em Vector Bundles on
Complex Projective Spaces}, Progress in Mathematics 3, Birkh\"{a}user, 1980

\bibitem[W]{W} Ward, R. S.: {\em On Self-Dual Gauge Fields}, Phys. Lett., {\bf 61} A,
      81-82 (1977)

\bibitem[WW]{WW} Ward, R. S. and Wells, R. O. : {\em Twistor Geometry and Field Theory},
      Cambridge University Press 1990

\end{thebibliography}
\end{document}